\documentclass[11pt]{article}
\usepackage[sectionbib]{natbib}
\usepackage{array,epsfig,fancyheadings,rotating}
\usepackage[colorlinks,linkcolor=blue,anchorcolor=blue,citecolor=blue,CJKbookmarks=True]{hyperref}
\usepackage{sectsty, secdot}
\usepackage{textcomp}
\usepackage{amsbsy}
\usepackage{amsmath}
\usepackage{amssymb}
\usepackage{epsfig}
\usepackage{multicol}
\usepackage{multirow}
\usepackage{color}
\usepackage{bm}
\usepackage{threeparttable}
\usepackage{float}
\usepackage{array}
\usepackage{amsmath}
\usepackage{amssymb}
\usepackage{amsfonts}
\usepackage{multirow}
\usepackage{amsthm}
\usepackage{graphicx}
\usepackage{bm}
\usepackage{enumerate}
\usepackage{subfigure}
\usepackage{booktabs}
\usepackage{pdflscape}
\usepackage{url}
\usepackage{xcolor,graphicx,float}
\usepackage{rotating}
\usepackage{times}
\usepackage{indentfirst}
\usepackage{blindtext}
\usepackage{algorithm,algorithmicx,float}
\usepackage{lipsum}
\usepackage[noend]{algpseudocode}
\usepackage{setspace}
\usepackage{booktabs}

\usepackage[sectionbib]{natbib}
\usepackage{array,epsfig,fancyheadings,rotating}
\usepackage[]{hyperref}  
\usepackage{sectsty, secdot}
\usepackage{authblk}
\usepackage{amsmath}
\usepackage{amssymb}
\usepackage{amsfonts}
\usepackage{multirow}
\usepackage{amsthm}
\usepackage{graphicx}
\usepackage{bm}
\usepackage{enumerate}
\usepackage{subfigure}
\usepackage{booktabs}
\usepackage{pdflscape}
\usepackage{url}
\usepackage{xcolor,graphicx,float}
\usepackage{rotating}
\usepackage{times}
\usepackage{indentfirst}
\usepackage{blindtext}
\usepackage{algorithm,algorithmicx,float}
\usepackage{lipsum}
\usepackage[noend]{algpseudocode}
\usepackage{setspace}
\usepackage{booktabs}


\makeatletter
\renewcommand{\@thesubfigure}{\hskip\subfiglabelskip}
\makeatother

\renewcommand{\baselinestretch} {1.5}
\makeatletter 
\def\singlespace{\def\baselinestretch{1}\@normalsize}

\newlength\savewidth

\@addtoreset{equation}{section}
\newtheorem{theorem}{Theorem}
\newtheorem{assumption}{Assumption}
\newtheorem{lemma}{Lemma}
\newtheorem{remark}{Remark}
\newtheorem{example}{Example}

\renewcommand{\hat}{\widehat}
\renewcommand{\tilde}{\widetilde}

\newcommand{\bg}{\begin{eqnarray}}
\newcommand{\ed}{\end{eqnarray}}
\newcommand{\bgn}{\begin{eqnarray*}}
\newcommand{\edn}{\end{eqnarray*}}

\makeatletter
\renewcommand{\@thesubfigure}{\hskip\subfiglabelskip}
\makeatother

\makeatletter
\def\singlespace{\def\baselinestretch{1}\@normalsize}

\date{\today}

\title{Minimax rates of convergence for nonparametric regression under adversarial attacks}

\author{Jingfu Peng}
\author{Yuhong Yang}

\affil{Yau Mathematical Sciences Center, Tsinghua University}

\begin{document}
\begin{sloppypar}

\maketitle

\begin{abstract}

Recent research shows the susceptibility of machine learning models to adversarial attacks, wherein minor but maliciously chosen perturbations of the input can significantly degrade model performance. In this paper, we theoretically analyse the limits of robustness against such adversarial attacks in a nonparametric regression setting, by examining the minimax rates of convergence in an adversarial sup-norm. Our work reveals that the minimax rate under adversarial attacks in the input is the same as sum of two terms: one represents the minimax rate in the standard setting without adversarial attacks, and the other reflects the maximum deviation of the true regression function value within the target function class when subjected to the input perturbations. The optimal rates under the adversarial setup can be achieved by an adversarial plug-in procedure constructed from a minimax optimal estimator in the corresponding standard setting. Two specific examples are given to illustrate the established minimax results.

\end{abstract}


\bigskip
\baselineskip=18pt

\section{Introduction}\label{sec:intro}

Over the past decade, machine/deep learning models have found unprecedented applications in a variety of domains including image recognition \citep{Krizhevsky2012}, natural language and speech processing \citep{Collobert2011}, game playing \citep{Silver2016}, autonomous driving \citep{Grigorescu2020}, many of which are safety-critical. However, it is found that these learning models are vulnerable to adversarial attacks. Here, an adversary is able to change the inputs to an already trained model, but cannot modify the training process. For example, input perturbations due to changes of weather conditions can significantly degrade the accuracy of trained neural networks for traffic sign recognition, demonstrating that such natural input variations present a significant challenge for deep learning \citep{Robey2020model}. Besides the nature as an adversary, a malicious opponent may choose perturbations to maximize prediction errors of a well trained neural network model \citep{Szegedy2014intriguing}. Similar vulnerabilities have been observed in various models across different application areas \citep[see, e.g.,][]{Biggio2013evasion, Goodfellow2015explaining, Papernot2016Limitations}.

The concerns about the safety and reliability of machine learning models have motivated a growing body of research focused on crafting the adversarial examples \citep{Goodfellow2015explaining, Papernot2016Limitations, Moosavi2016DeepFool, Carlini2017Towards, Awasthi2020adversarial} and devising defenses to enhance model robustness against such perturbations \citep{Goodfellow2015explaining, Madry2018towards, FINLAY2021100017, Raghunathan2018certified, Cohen2019}. Adversarial training, which minimizes the empirical risk under worst-case perturbations on the training data, has been empirically demonstrated to be effective against various attacks \citep[see, e.g.,][]{Madry2018towards}. While considerable efforts have been made on constructing attack and defence, the problem of understanding the intrinsic hardness of estimation and assessing the optimality of learning methods under adversarial attacks are far less understood.

One of the most important approach to measuring the difficulty of a nonparametric statistical problem is to evaluate its minimax risk \citep[see, e.g.,][]{ibragimov1982bounds, birge1986estimating, Yang1999information}. In the adversarial setting, the maximal risk of an estimator is defined as its worst statistical performance over a class of distributions when the input perturbation is generated from a given perturbation set to deprave the model's performance. If its maximal risk is minimal (rate) among all estimators, then this estimator is called minimax (rate) optimal. To the best of our knowledge, investigating the adversarial robustness from the minimax viewpoint has not been paid much attention. \cite{Dan2020sharp} considered a binary classification problem with data generated from a Gaussian mixture model. They established the minimax rate of excess risk when the perturbations lie in an origin-symmetric convex set. \cite{Xing2021predictive} determined the minimax rate of a nonparametric classification problem when the testing input is randomly perturbed on a sphere, and established the minimax optimality of a nearest neighbor rule. In a setup of linear regression with Gaussian regressors, \cite{Xing2021adversarially} provided the minimax rate for estimating regression coefficients under bounded $\ell_2$-norm perturbations. In a context of data contamination where a subset of training sample can be arbitrarily modified by an attacker, \cite{zhao2023robust} established the minimax rates for the estimation of a nonparametric Lipschitz regression function under both $\ell_2$ and $\ell_\infty$ losses. Although the above theoretical advancements provide valuable insights, they are confined to some restricted setups based on simple models and architectures, and thus do not seem to be applicable to the broader nonparametric setting with the adversarial attacks as we consider.

Under a nonparametric regression setting with minimal assumptions regarding the adversarial perturbations, an important question arises: What is the minimax rate of convergence for a general class of regression functions?

This paper determines the sup-norm rate of convergence in a nonparametric regression setup with additive perturbations, in which the attacker can add arbitrary perturbations in a set to the input, thereby degrading the performance of the trained estimator. We establish that under general class of regression functions and adversarial perturbation sets, the minimax risk converges at the order of the rate in the standard setup without adversaries, plus the maximum deviation of true function values within the target regression function class. The optimal rate can be achieved by an adversarial plug-in procedure constructed from a minimax optimal estimator in the standard setting. We provide minimax results for two specific examples of function classes, including isotropic H\"{o}lder class and anisotropic H\"{o}lder class, and investigate the effects of $\ell_p$-attacks ($0< p \leq \infty$) and sparse attacks under these two function classes, respectively.

\subsection{Related work}

\emph{Sup-norm convergence.} Determining the rate of convergence in the sup-norm is a crucial topic in statistics and machine learning. Classical contributions include works by \cite{Devroye1978, Stone1982optimal, donoho1994asymptotic, Korostelev1999, lepski2000asymptotically, bertin2004asymptotically, GAIFFAS2007782, Gine2009, CHEN2015447}. More recently, the implications of sup-norm convergence in transfer learning have been explored by \cite{schmidt2022local}, and its relation to adversarial training has been investigated by \cite{imaizumi2023sup}. However, these studies focus on standard setups without adversarial perturbations to the input data.

\emph{Robustness of nonparametric classifiers.} Several previous works analyzed the robustness of specific families of classifiers. \cite{Wang2018Analyzing} studied the robustness of nearest neighbor classifier. \cite{Yang2020Robustness} proposed the attack strategies that apply to a wide range of non-parametric classifiers and analyzed a general defense method based on data pruning. \cite{bhattacharjee2020non} proved the consistency of the nearest neighbor and kernel estimators. Note that the aforementioned works do not establish the optimal rate of convergence of nonparametric estimation under the adversarial attacks.

\emph{Distributional robustness optimization. }\cite{Lee2018Minimax} and \cite{Tu2019} established the connections between the adversarial training and distributional robustness optimization (DRO) \citep{ben2009robust, shapiro2021lectures}. These connections can be used to upper bound the generalization error of the adversarial training. In the context of DRO, when the loss function is defined as a product of the response variable and the parameter, \cite{Duchi2023Latent} obtained minimax lower bounds for a distributionally robust loss. However, the linear form of the loss function in their work cannot be applied to the typical regression setting.

\emph{Other related work. }Rather than studying the minimax risk, another line of work obtained tight statistical characterizations of the Bayes adversarial risk and developed classifiers to realized it \citep{Schmidt2018adversarially, Bhagoji2019lower, Pydi2020adversarial}. The trade-offs between standard and robust accuracy have been studied by \cite{Madry2018towards, Schmidt2018adversarially, Tsipras2019robustness, Raghunathan2019adversarial, Zhang2019theoretically, Javanmard2020precise, Min2021curious, Mohammad2021, Dobriban2020provable, Javanmard2022precise}. Algorithm-free generalization bounds such as VC-dimension have been studied by \cite{Attias2019improved, Montasser2019vc} in the adversarial setting. Rademacher complexity of the adversarial training has been investigated by \cite{Yin2019rademacher, khim2018adversarial, Awasthi2020adversarial}. Recently, \cite{liu2023non} derived non-asymptotic bounds for adversarial excess risk under misspecified models. Note that the above analyses primarily center on upper bounding the adversarial risk, thus lacking corresponding lower bounds necessary for determining the minimax rates.

\subsection{Outline}

The rest of this paper is organized as follows. Section~\ref{sec:setup} gives a setup for the nonparametric regression problem and the definition of adversarial loss/risk. In Section~\ref{sec:main}, we state upper and lower bounds on the minimax risks under the adversarial attack. Two specific examples are discussed in Section~\ref{sec:example}. Section~\ref{sec:simulation} presents numerical simulation results. The proofs of the main theorems and examples are provided in the Appendix.

\section{Problem setup}\label{sec:setup}

This paper considers the problem of nonparametric regression estimation. Suppose the observations $(X_1,Y_1),\ldots,(X_n,Y_n) \in \mathcal{X}\times \mathcal{Y}$ are generated from the regression model
\begin{equation}\label{eq:model}
  Y_i = f(X_i) + \xi_i,
\end{equation}
where $\mathcal{X}\subseteq \mathbb{R}^d$, $\mathcal{Y} \subseteq \mathbb{R}$, $f:\mathcal{X}\to \mathcal{Y}$ is an unknown regression function, $\xi_i$ is a random error term with $\mathbb{E}(\xi_i| X_i)=0$ a.s., and $X_i$ follows an unknown marginal distribution $\mathbb{P}_{X}$ on $\mathcal{X}$. The goal is to develop an estimator $\hat{f}$ of $f$ based on the observed data. The estimation accuracy of $\hat{f}$ is measured by the sup-norm loss. In the standard setting of regression with unperturbed future $X$ values, this loss is defined as $\sup_{x \in \mathcal{X}}| f(x) - \hat{f}(x) |$, which  quantifies the uniform convergence of $\hat{f}$ to $f$ over $\mathcal{X}$.

In this paper, we consider the estimation of the regression function in the presence of an adversary. Specifically, when assessing the performance of the estimator $\hat{f}$, the adversary can add any perturbation $\delta \in \Delta_n$ to the input $x$, where $\Delta_n \in \mathbb{R}^d$ is a closed set containing $ \delta = 0$, and $\Delta_n$ may depend on the sample size $n$. A representative example of $\Delta_n$ is the $\ell_p$-ball $B_p^{q_n}=\{ z: \| z \|_p \leq q_n \}$ centering at origin with radius $q_n>0$ and $p>0$. In the adversarial setting, the sup-norm loss of estimation is defined as
\begin{equation}\label{eq:adversarial_loss}
  L_{\Delta_n}(f,\hat{f})=\sup_{x \in \mathcal{X}}\sup_{\substack{\delta\in \Delta_n\\x+\delta \in \mathcal{X}}}\left|f(x)-\hat{f}(x+\delta)\right|,
\end{equation}
and the corresponding adversarial risk is given by
\begin{equation}\label{eq:adversarial_risk}
  R_{\Delta_n}(f,\hat{f})=\mathbb{E} L_{\Delta_n}(f,\hat{f}),
\end{equation}
where the expectation $\mathbb{E}$ is taken with respect to the observed data generated from the regression model (\ref{eq:model}), and the subscript $\Delta_n$ here is employed to emphasize the dependence of the adversarial risk/loss on the perturbation set $\Delta_n$. In the standard regression setting with $\Delta_n = \{ 0\}$, expressions (\ref{eq:adversarial_loss}) and (\ref{eq:adversarial_risk}) reduce to the standard sup-norm loss
$$
L(f,\hat{f}) = \sup_{x \in \mathcal{X}}\left| f(x) - \hat{f}(x) \right|
$$
and the standard sup-norm risk
$$
R(f,\hat{f})=\mathbb{E}L(f,\hat{f}),
$$
respectively. In the adversarial setting, an estimator $\hat{f}$ is sought to be robust to the adversarial perturbation of $x$.

The regression function $f$ is assumed to belong to a function class $\mathcal{F}$. The minimax risk of estimating $f \in \mathcal{F}$ under the adversarial sup-norm loss is expressed as:
\begin{equation}\label{eq:minimax_V}
  V_{\Delta_n}=\inf_{\hat{f}}\sup_{f \in\mathcal{F}}R_{\Delta_n}(f,\hat{f}).
\end{equation}
Then two important questions arise:
\begin{description}
  \item[Q1.] What factors determine the rate of convergence of $V_{\Delta_n}$?
  \item[Q2.] How can minimax optimal procedures be developed to achieve the optimal rate of $V_{\Delta_n}$?
\end{description}
Answers to questions Q1 and Q2 have the potential to offer previously unavailable insights into the theoretical foundations and practical applications of adversarial learning.

Throughout this paper, let $\mathbb{N}_0$ denote the set of non-negative integers. For any $a \in \mathbb{R}^d$ and $B \subseteq \mathbb{R}^d$, we use the Minkowski sum notations $a+B\triangleq \{ a+b:b\in B\}$ and $a-B\triangleq \{ a-b:b\in B\}$. For any positive sequences $a_n$ and $b_n$, we denote $a_n = O(b_n)$ and $a_n \lesssim b_n$ if there exist $C>0$ and $N>0$ such that $n\geq N$ implies $a_n \leq Cb_n$. If $a_n= O(b_n)$ and $b_n=O(a_n)$, then we write $a_n\asymp b_n$. For $1 \leq p < \infty$, we use $\|\delta\|_p$ to denote the $\ell_p$-norm $(\sum_{j=1}^{d}|\delta|_j^p)^{1/p}$ of a vector $\delta \in \mathbb{R}^d$. We use $\|\delta\|_\infty$ to denote the sup-norm $\sup_{1 \leq j \leq d}|\delta_j|$. For brevity, we write $\|\delta\|$ to represent the $\ell_2$-norm.

\section{Main results}\label{sec:main}

In this section, we begin by deriving a closed form expression for the ideal adversarial loss $\inf_{f'}L_{\Delta_n}(f,f')$. Then we establish the minimax rates of convergence for the general function classes $\mathcal{F}$ and perturbation sets $\Delta_n$.

\subsection{Ideal adversarial loss}\label{sec:sub:optimal}

We first introduce an equivalent form for the adversarial sup-norm loss (\ref{eq:adversarial_loss}), which offers conveniences in characterizing both the ideal adversarial loss and the minimax risk $V_{\Delta_n}$.

\begin{lemma}\label{lem:equ_form}
  For any estimator $\hat{f}$, we have
  \begin{equation}\label{eq:equ_form}
    L_{\Delta_n}(f,\hat{f}) = \sup_{x \in \mathcal{X}}\sup_{x'\in (x+\Delta_n)\cap \mathcal{X}}\left|f(x)-\hat{f}(x')\right| = \sup_{x' \in \mathcal{X}}\sup_{x \in (x'-\Delta_n)\cap \mathcal{X}} \left| f(x) - \hat{f}(x') \right|.
  \end{equation}
\end{lemma}

Lemma~\ref{lem:equ_form} provides an alternative expression for the adversarial loss by exchanging the order of two supremum operations. The inner supremum in the last argument of (\ref{eq:equ_form}), which depends on the perturbation set, is taken respect to the regression function $f$ rather than the estimator $\hat{f}$. This property facilitates the derivation of the ideal adversarial loss and the ideal adversarial estimator (i.e., the best performing ``estimate'' when the underlying regression function $f$ is known). The next theorem addresses this aspect.

\begin{theorem}\label{theo:oracle}
  Given the regression function $f$, the ideal adversarial loss is given by
  \begin{equation}\label{eq:oracle_loss}
    L_{\Delta_n}^*(f) \triangleq \inf_{f'}L_{\Delta_n}(f,f')=\frac{1}{2}\sup_{x' \in \mathcal{X}}\left[ \sup_{x \in (x'-\Delta_n)\cap \mathcal{X}}f(x)-\inf_{x \in (x'-\Delta_n)\cap \mathcal{X}}f(x)\right],
  \end{equation}
  where the minimum is achieved by the adversarial regression function:
  \begin{equation}\label{eq:oracle_estimator}
    f^*(x)= \frac{1}{2}\left[ \sup_{x' \in (x-\Delta_n)\cap \mathcal{X}}f(x')+\inf_{x' \in (x-\Delta_n)\cap \mathcal{X}}f(x')\right], \quad x \in \mathcal{X}.
  \end{equation}
\end{theorem}

Theorem~\ref{theo:oracle} provides a closed form expression for the ideal adversarial loss, which shows that the ideal adversarial loss is proportional to the maximum variation of the true regression function value within the perturbation set $\Delta_n$ over the domain $\mathcal{X}$. Moreover, the ideal adversarial regression function is exactly the average of the maximum and minimum values of the function $f$ in the adversarial neighborhood $(x-\Delta_n)\cap \mathcal{X}$.

The result from Theorem~\ref{theo:oracle} substantiates that the optimal adversarial robustness is jointly determined by the size of the perturbation set and the smoothness of the true regression function. For example, when $f$ satisfies the Lipschitz smoothness condition $|f(x) - f(z)| \leq L\cdot\|x - z\|$ and $\Delta_n$ has the diameter $\mathrm{diag}(\Delta_n)\triangleq \max_{x,z}\|x-z \|$, then the ideal adversarial loss $$
L_{\Delta_n}^*(f) \leq  \frac{L\cdot\mathrm{diag}(\Delta_n)}{2},
$$
a quantity controllable when the diameter of $\Delta_n$ is not excessively large. In contrast, if the true regression function is discontinuous, then $L_{\Delta_n}^*(f)$ cannot degenerate to $0$ unless $\Delta_n=\{0 \}$. Also, if $\Delta_n$ does not shrink with $n$, $L_{\Delta_n}^*(f)$ may not converge to $0$.

\begin{remark}
  In the literature, several papers have obtained precise characterizations or tight bounds on the ideal adversarial loss \citep[see, e.g.,][]{Bhagoji2019lower, Pydi2020adversarial, Dan2020sharp, Xing2021adversarially}. However, it is important to note that all of these works focus on parametric models, which cannot imply the adversarial robustness for nonparametric regression as considered in this paper.
\end{remark}

\subsection{Minimax rates of convergence}\label{sec:sub:minimax}

In this subsection, our aim is to establish the minimax rates of convergence for the sup-norm risk under the adversarial attacks. We propose an adversarial plug-in procedure to achieve the minimax optimal rates, which is derived from a minimax optimal estimator in the corresponding standard setting.

In Theorem~\ref{theo:oracle}, we obtain the explicit expression for the ideal adversarial regression function (\ref{eq:oracle_estimator}). However, (\ref{eq:oracle_estimator}) is infeasible in practice as it relies on the true regression function $f$. Motivated by (\ref{eq:oracle_estimator}), we devise a feasible adversarial estimator through the following two steps:
\begin{description}
  \item[Step 1.] Utilizing the observed data $(X_1,Y_1),\ldots,(X_n,Y_n)$, we construct an estimator $\tilde{f}$ for the regression function $f$.
  \item[Step 2.] Subsequently, we formulate an adversarial plug-in estimator:
  \begin{equation}\label{eq:plug_in}
  \hat{f}_{\mathrm{PI}}(x) = \frac{1}{2}\left[ \sup_{x' \in (x-\Delta_n)\cap \mathcal{X}}\tilde{f}(x')+\inf_{x' \in (x-\Delta_n)\cap \mathcal{X}}\tilde{f}(x')\right], \quad x \in \mathcal{X}.
  \end{equation}
\end{description}

The performance of the adversarial plug-in estimator $\hat{f}_{\mathrm{PI}}(x)$ clearly depends on the construction of $\tilde{f}$. The following theorem first provides an upper bound for the adversarial risk of $\hat{f}_{\mathrm{PI}}(x)$ considering a general $\tilde{f}$. Additionally, Theorem~\ref{theo:upper} establishes minimax upper bounds when specific choices of $\tilde{f}$ are adopted.

\begin{theorem}[Upper bound]\label{theo:upper}
  For any regression function $f$ and any estimator $\tilde{f}$, the adversarial risk of the plug-in estimator (\ref{eq:plug_in}) is upper bounded by
  \begin{equation}\label{eq:upper}
    R_{\Delta_n}(f,\hat{f}_{\mathrm{PI}}) \leq R(f,\tilde{f}) + L_{\Delta_n}^*(f),
  \end{equation}
  where $L_{\Delta_n}^*(f)$ is the ideal adversarial loss defined in (\ref{eq:oracle_loss}).

  Moreover, given a function class $\mathcal{F}$, if $\tilde{f}$ satisfies
  \begin{equation}\label{eq:minimax_stand}
    \sup_{f \in \mathcal{F}}R(f,\tilde{f}) \asymp \inf_{\hat{f}}\sup_{f \in \mathcal{F}}R(f,\hat{f}),
  \end{equation}
  then the adversarial maximal risk of $\hat{f}_{\mathrm{PI}}$ is upper bounded by
  \begin{equation}\label{eq:minimax_upper}
    \sup_{f \in \mathcal{F}}R_{\Delta_n}(f,\hat{f}_{\mathrm{PI}}) \lesssim \inf_{\hat{f}}\sup_{f \in \mathcal{F}}R(f,\hat{f}) +  \sup_{f \in \mathcal{F}}L_{\Delta_n}^*(f).
  \end{equation}

\end{theorem}

The relationship (\ref{eq:upper}) illustrates that the adversarial risk of the plug-in estimator $\hat{f}_{\mathrm{PI}}$ can be upper bounded by the standard risk of the original estimator $\tilde{f}$ plus a multiple of the ideal adversarial loss $L_{\Delta_n}^*(f)$. Importantly, this relation holds without any additional constraints on the true regression function and the perturbation set, and without imposing assumptions on the estimator $\tilde{f}$. The second part of Theorem~\ref{theo:upper} indicates that if the original estimator $\tilde{f}$ is minimax optimal in the standard setting, then the corresponding adversarial maximal risk $\sup_{f \in \mathcal{F}}R_{\Delta_n}(f,\hat{f}_{\mathrm{PI}})$ is upper bounded by the standard minimax rate plus $\sup_{f \in \mathcal{F}}L_{\Delta_n}^*(f)$.

The following lower bound results show that the adversarial plug-in estimator based on $\tilde{f}$ with (\ref{eq:minimax_stand}) is in fact minimax rate optimal.
\begin{theorem}[Lower bound]\label{theo:lower}
  For any regression function $f$ and any estimator $\hat{f}$, the adversarial risk is lower bounded by
  \begin{equation}\label{eq:lower}
    R_{\Delta_n}(f,\hat{f}) \geq R(f,\hat{f})\vee L_{\Delta_n}^*(f).
  \end{equation}

  Furthermore, for any function class $\mathcal{F}$, we have
  \begin{equation}\label{eq:minimax_lower}
    \inf_{\hat{f}}\sup_{f \in \mathcal{F}}R_{\Delta_n}(f,\hat{f}) \gtrsim \inf_{\hat{f}}\sup_{f \in \mathcal{F}}R(f,\hat{f})+ \sup_{f \in \mathcal{F}} L_{\Delta_n}^*(f).
  \end{equation}
\end{theorem}

In summary, Theorems~\ref{theo:upper}--\ref{theo:lower} together establish the minimax rates of convergence for nonparametric regression under the adversarial attacks,
\begin{equation}\label{eq:minimax}
    \inf_{\hat{f}}\sup_{f \in \mathcal{F}}R_{\Delta_n}(f,\hat{f}) \asymp \inf_{\hat{f}}\sup_{f \in \mathcal{F}} R(f,\hat{f})+  \sup_{f \in \mathcal{F}}L_{\Delta_n}^*(f).
\end{equation}
Therefore, (\ref{eq:minimax}) addresses Question Q.1 raised in Section~\ref{sec:setup}, showing that the adversarial minimax rate is jointly determined by the standard minimax rate and the largest ideal loss in $\mathcal{F}$. Regarding Question Q.2, we establish that if $\tilde{f}$ is minimax optimal in the sense that $\sup_{f \in \mathcal{F}}R(f,\tilde{f}) \asymp \inf_{\hat{f}}\sup_{f \in \mathcal{F}} R(f,\hat{f})$ under the standard setting, then the adversarial plug-in estimator $\hat{f}_{\mathrm{PI}}$ based on $\tilde{f}$ is minimax optimal in terms of the adversarial risk. To the best our knowledge, (\ref{eq:minimax}) is the first minimax result in adversarial learning for the general regression setting. Our bounds are modular and can be applied to many models by computing the sup-norm convergence and the ideal adversarial loss in the target function class.

\section{Applications}\label{sec:example}

In this section, we demonstrate the applications of the theorems in the previous section through specific examples of function classes and perturbation sets. We consider the case $\mathcal{X}=[0,1]^d$, and $(X_1, Y_1), \ldots, (X_n, Y_n)$ are drawn i.i.d.\ according to the regression model (\ref{eq:model}). The following assumption on the distribution of $X$ is required.

\begin{assumption}\label{ass:density}
  The marginal distribution $\mathbb{P}_{X}$ admits a density function that is lower bounded away from 0 and upper bounded by a positive constant on $\mathcal{X}$.
\end{assumption}

Assumption~\ref{ass:density} ensures that the covariates $X$ are more or less evenly distributed over the compact support $[0,1]^d$. As a result, there are sufficiently many observations around any point in the support, allowing for the construction of well-behaved estimators for the regression function in the sup-norm loss. This assumption is standard in nonparametric regression with random design; see, for example, Condition~3' in \cite{Stone1982optimal} and Definition~2.2 in \cite{Audibert2007Fast}. In addition, we further assume that the random error term is distributed according to a centered Gaussian distribution, which is the scenario where the known minimax theory in sup-norm can apply \citep[see, e.g.,][]{Stone1982optimal, Bertin2004minimax, GAIFFAS2007782}.

\begin{assumption}\label{ass:rand_error}
  The random error term $\xi$ follows a zero-mean Gaussian distribution and is independent of $X$.
\end{assumption}

\subsection{Isotropic H\"{o}lder class}\label{sec:Isotropic}

Let $\beta=k+\alpha$ for some $k \in \mathbb{N}_0$ and $0<\alpha\leq 1$, and let $L>0$. A function $f:[0,1]^d \to \mathbb{R}$ called $(\beta,L)$-smooth if for every $(k_1,\ldots,k_d)$, $k_i \in \mathbb{N}_0$, and $\sum_{i=1}^{d}k_i=k$, the partial derivative $\partial^k f/(\partial x_1^{k_1}\cdots \partial x_d^{k_d})$ exists and satisfies
\begin{equation}\label{eq:iso_holder}
  \left| \frac{\partial^k f}{\partial x_1^{k_1}\cdots \partial x_d^{k_d}}(x) - \frac{\partial^k f}{\partial x_1^{k_1}\cdots \partial x_d^{k_d}}(z) \right|\leq L \cdot \left\| x - z \right\|^{\alpha}
\end{equation}
for all $x,z \in [0,1]^d$. The isotropic H\"{o}lder class, denoted $\mathcal{F}_1(\beta,L)$, is defined as the set of all $(\beta,L)$-smooth functions $f:[0,1]^d \to \mathbb{R}$.

\begin{example}\label{exam:iso}
Suppose Assumptions~\ref{ass:density}--\ref{ass:rand_error} are satisfied. For any closed perturbation set $\Delta_n \in \mathbb{R}^d$, define
\begin{equation}\label{eq:r_iso}
r_n \triangleq \max_{\delta_1, \delta_2 \in \Delta_n }\|\delta_1 - \delta_2 \|.
\end{equation}
If there exists a pair of $\delta$ and $\delta'$ in $\Delta_n$ such that $\|\delta - \delta' \|=r_n$ and $\{t\delta  + (1-t)\delta':0\leq t \leq 1 \} \subseteq \Delta_n$, then we have
\begin{equation}\label{eq:rate_iso}
    \inf_{\hat{f}}\sup_{f \in \mathcal{F}_1(\beta,L)}R_{\Delta_n}(f,\hat{f}) \asymp \left(\frac{\log n}{n}\right)^{\frac{\beta}{2\beta+d}} + C_d r_n^{ 1 \wedge\beta},
\end{equation}
where $C_d \leq C d^{k/2}$ is a constant not depending on $n$.
\end{example}

In view of (\ref{eq:minimax}), the proof of the result in Example~\ref{exam:iso} consists of examining the standard minimax rate $\inf_{\hat{f}}\sup_{f \in \mathcal{F}_1(\beta,L)} R(f,\hat{f})$ and the rate of $\sup_{f \in \mathcal{F}_1(\beta,L)}L_{\Delta_n}^*(f)$. The standard minimax rate within the isotropic H\"{o}lder class is established in \cite{Stone1982optimal}, which demonstrates that
\begin{equation}\label{eq:rate_iso_stan}
  \inf_{\hat{f}}\sup_{f \in \mathcal{F}_1(\beta,L)} R(f,\hat{f}) \asymp \left(\frac{\log n}{n}\right)^{\frac{\beta}{2\beta+d}}.
\end{equation}
The determination of the rate of $\sup_{f \in \mathcal{F}_1(\beta,L)}L_{\Delta_n}^*(f)$ is provided in Section~S1 of the Supplementary Material.

The quantity $r_n$ in (\ref{eq:r_iso}) measures the length of the longest line segment contained in the set $\Delta_n$, and it may depend on the sample size $n$. The condition imposed on $\Delta_n$ is quite mild, which is satisfied by the $\ell_p$-ball: $B_p^{q_n} \triangleq \{ \delta\in \mathbb{R}^d: \| \delta \|_p \leq q_n \}$, $0<p\leq\infty$, and the $\ell_p$-ball with the $\ell_0$-constraint: $B_p^{q_n}\cap \{\delta: \| \delta \|_0 \leq s_n \}$. Note that there is an extensive body of prior work studying adversarial machine learning based on $\ell_0$ \citep{Delgosha2024}, $\ell_2$ \citep{bhattacharjee2020non, bhattacharjee2021sample}, and $\ell_\infty$ attacks \citep{Athalye2018Obfuscated, marzi2018sparsity}. However, these analyses focus on the specific attacks and lack general applicability. In contrast, the result in Example~\ref{exam:iso} sheds theoretical insight on the adversarial robustness under the general $\ell_p$-attacks with $0< p \leq \infty$. Specifically, when $\Delta_n = B_p^{q_n}$, we have $r_n = q_n$, and thus the minimax adversarial risk is given by
\begin{equation}\label{eq:rate_iso_ball}
  \inf_{\hat{f}}\sup_{f \in \mathcal{F}_1(\beta,L)}R_{B_p^{q_n}}(f,\hat{f}) \asymp \left(\frac{\log n}{n}\right)^{\frac{\beta}{2\beta+d}} + C_d{q_n}^{ 1 \wedge\beta},
\end{equation}
which can be reached by the adversarial plug-in estimator (\ref{eq:plug_in}) with $\tilde{f}$ constructed by a suitably designed local polynomial estimator (see, e.g., \cite{Stone1982optimal}, \cite{GAIFFAS2007782}, and \cite{tsybakov2008introduction}).

The equation~(\ref{eq:rate_iso_ball}) shows that when $\beta < 1$ and $q_n \lesssim (\log n/ n)^{1/(2\beta+d)}$, the minimax rate in the adversarial sup-norm remains unchanged to the standard minimax rate (\ref{eq:rate_iso_stan}). However, as the magnitude of perturbation increases, e.g., $q_n \gtrsim (\log n/ n)^{1/(2\beta+d)}$, the minimax risk has the order $q_n^\beta$. When $\beta \geq 1$ and the functions in $\mathcal{F}_1(\beta,L)$ become smoother, the critical radius $q_n$ for the phase transition is $(\log n/ n)^{\beta/(2\beta+d)}$. It is also worth noting that the norm parameter $p$, which controls the shape of the perturbation set $B_p^{q_n}$, does not affect the adversarial minimax rates in this example. This is because the adversarial risk $R_{B_p^{q_n}}(f,\hat{f})$ is defined via the worst-case perturbation within the ball $B_p^{q_n}$, and the minimax risk considers the worst-case adversarial risk over all functions in the function class. In fact, the maximum adversarial risk over the function class $\mathcal{F}_1(\beta,L)$ is attained at certain functions evaluated at points $x',x$ satisfying $\|x' - x\|_2 = q_n$. Therefore, in the minimax sense, $p$ does not influence the adversarial minimax risk. However, in other regression function classes of interest, the shape of the perturbation may have an effect on the robustness of a given estimator; see Section~\ref{sec:anisotropic} for further discussion.

\begin{remark}
  In this paper, we primarily focus on the adversarial sup-norm as the robustness performance measure. Using the uniformity of the sup-norm loss, we can derive the following upper bound on the adversarial $L_2$-loss
  \begin{equation*}
    \bar{L}_{\Delta_n}(f,\hat{f}) \triangleq \int_{\mathcal{X}}\sup_{\substack{\delta\in \Delta_n\\x+\delta \in \mathcal{X}}} \left| f(x) - \hat{f}(x+\delta) \right|^2 \mathbb{P}_{X}dx \lesssim L_{\Delta_n}^2(f,\hat{f}),
  \end{equation*}
  under the assumption that $\mathcal{X}$ is a compact set and $\mathbb{P}_X$ satisfies Assumption~\ref{ass:density}. Based on this relation and (\ref{eq:rate_iso}), we can also derive an upper bound on the minimax adversarial risk under $L_2$-loss over the isotropic H\"{o}lder class:
  \begin{equation*}
    \left(\frac{\log n}{n}\right)^{\frac{2\beta}{2\beta+d}} + C_d{r_n}^{ 2(1 \wedge\beta)}.
  \end{equation*}
  It remains to be seen if this is the minimax optimal rate.
\end{remark}

\subsection{Anisotropic H\"{o}lder class}\label{sec:anisotropic}

In practice, one of the typically desired properties of a regression function or its estimator is that it is invariant or robust against changes or perturbations of an input in some specific directions. For example, in image classification tasks, the target function should be invariant against a spatial shift or rotation of an input image \citep{Simard2003, Krizhevsky2012}. In the same spirit, in the context of autonomous driving, a traffic sign recognition model should be trained to be robust to natural variations in severe weather conditions.

Motivated by these examples, in this subsection, we investigate the adversarial minimax risks on the anisotropic H\"{o}lder class $\mathcal{F}_2(\beta,L)$, where $\beta=(\beta_1,\ldots,\beta_d)\in (0,1]^d$ and $L=(L_1,\ldots,L_d) \in (0,\infty)^d$ \citep{birge1986estimating, bertin2004asymptotically, bhattacharya2014anisotropic, jeong2023art}. This class is defined by
\begin{equation}\label{eq:aniso_holder}
\begin{split}
      \mathcal{F}_2(\beta,L) \triangleq \left\{ f: [0,1]^d \to \mathbb{R}:\right. & \left| f(x) - f(z) \right| \\
     & \left. \leq L_1|x_1-z_1|^{\beta_1}+\cdots+ L_d|x_d-z_d|^{\beta_d} \right\},
\end{split}
\end{equation}
which is a set of functions that have ``direction-dependent'' smoothness, whereas the isotropic H\"{o}lder class considered in Section~\ref{sec:Isotropic} assumes isotropic smoothness that is uniform in all directions.

\begin{example}\label{exam:aniso}
  Suppose Assumptions~\ref{ass:density}--\ref{ass:rand_error} hold. For any perturbation set $\Delta_n \in \mathbb{R}^d$, define $r_i \triangleq \sup_{\delta,\delta' \in \Delta_n}|\delta_i - \delta_i'|$ for $1\leq i \leq d$, where $\delta=(\delta_1,\ldots,\delta_d)$ and $\delta'=(\delta'_1,\ldots,\delta'_d)$. Then we have
  \begin{equation}\label{eq:rate_aniso}
    \inf_{\hat{f}}\sup_{f \in \mathcal{F}_2(\beta,L)}R_{\Delta_n}(f,\hat{f}) \asymp \left(\frac{\log n}{n}\right)^{\frac{\bar{\beta}}{2\bar{\beta}+d}} + \max\left\{r_1^{\beta_1},\ldots, r_d^{\beta_d}\right\},
  \end{equation}
  where $\bar{\beta}=d/(\sum_{i=1}^{d}1/\beta_i)$.
\end{example}

The first term on the right side of (\ref{eq:rate_aniso}) represents the standard minimax rate under the sup-norm, which is determined by the average smoothness and the dimension $d$. The second term is related to the maximum deviation of function values along each coordinates. Combining the results in Section~\ref{sec:main} with \cite{bertin2004asymptotically, Bertin2004minimax}, it can be deduced that the adversarial minimax rate is achievable through the plug-in estimator (\ref{eq:plug_in}), with $\tilde{f}$ being a multivairate kernel estimator with different bandwidths across different coordinates.

To compare the adversarial minimax rates in the isotropic and anisotropic H\"{o}lder classes, let us consider a specific perturbation set $\Delta_n=\{\delta: |\delta_1| \leq q_n,\delta_2=\cdots=\delta_d=0 \}$, where $q_n\to 0$ and $q_n\gtrsim (\log n/n)^{1/(2\bar{\beta}+d)}$. Note that the attacks within $\Delta_n$ are concentrated solely on the first coordinate. Suppose $\beta_1 > \bar{\beta}$. The isotropic H\"{o}lder class with the smoothness parameter $\bar{\beta}$ exhibits the minimax rate:
$$
\inf_{\hat{f}}\sup_{f \in \mathcal{F}_1(\bar{\beta},L)}R_{\Delta_n}(f,\hat{f})\asymp q_n^{\bar{\beta}}.
$$
In contrast, for the anisotropic H\"{o}lder class, the minimax rate is:
$$
\inf_{\hat{f}}\sup_{f \in \mathcal{F}_2(\beta,L)}R_{\Delta_n}(f,\hat{f})\asymp \max\{r_1^{\beta_1},\ldots, r_d^{\beta_d}\} = q_n^{\beta_1},
$$
which converges significantly faster than $\inf_{\hat{f}}\sup_{f \in \mathcal{F}_1(\bar{\beta},L)}R_{\Delta_n}(f,\hat{f})$ as $q_n^{\beta_1}/q_n^{\bar{\beta}} \to 0$. This phenomenon implies that although the average smoothness is the same for the two function classes, when the attack is only in a smoother direction, the adversarial minimax risk in the anisotropic H\"{o}lder class is faster than that in the isotropic H\"{o}lder class.

\section{Simulation studies}\label{sec:simulation}

In this section, we present several numerical experiments to illustrate the theoretical results established in Sections~\ref{sec:main}--\ref{sec:example}. The data are generated from the model (\ref{eq:model}), where $\mathcal{X}=[0,1]^2$, $X$ follows a uniform distribution on $[0,1]^2$, and $\xi$ is independent of $X$ and distributed as $N(0,\sigma^2)$. We consider several regression functions and attack scenarios:
\begin{description}
  \item[Case 1] $f(x_1,x_2) = \sqrt{x_1x_2}$ with perturbation set $\Delta_n = B_{\infty}^r$.
  \item[Case 2] $f(x_1,x_2) = \sqrt{(x_1 - 0.5)^2+(x_2 - 0.5)^2}$, with $\Delta_n = B_{\infty}^r$.
  \item[Case 3] $f(x_1,x_2) = \sqrt{x_1} + 0.1x_2 - 0.5$, with $\Delta_n = [-4r, + 4r] \times [-r/4, + r/4]$.
  \item[Case 4] $f(x_1,x_2) = \sqrt{x_1} + 0.1x_2 - 0.5$, with $\Delta_n = [-r/4, + r/4] \times [-4r, + 4r]$.
\end{description}
In each case, $\sigma^2$ is adjusted so that the signal-to-noise ratio equals $5$. The attack magnitude $r$ increases from 0 to 0.1. Cases 1--2 serve as two representative examples of isotropic H\"{o}lder classes, where the perturbation set is chosen as the $\ell_\infty$-ball. In contrast, Cases 3--4 consider regression functions with different degrees of variation along different axes, where the attack magnitudes are also anisotropic.

We consider three competing methods. The baseline method (LP) is the classical local polynomial regression studied in \cite{Stone1982optimal}, \cite{Bertin2004minimax}, and \cite{GAIFFAS2007782} based on the rectangular kernel. We employ a polynomial of degree $\ell = 1$ (i.e., local linear regression). In Cases 1 and 2, the bandwidth is set as $h = n ^{-1/(0.5+2)}$ and $h = n^{-1/(1+2)}$, respectively. In Cases 3 and 4, we use different bandwidths for different coordinates, setting $h_1 = n ^{-1/(0.5+2)}$ and $h_2 = n^{-1/(1+2)}$. These choices are theoretically proven to achieve the standard minimax rates in the respective cases.

\begin{figure}[!t]
    \centering
    \subfigure[Case 1]{
    \begin{minipage}[t]{1\linewidth}
    \centering
       \includegraphics[width=5.3in]{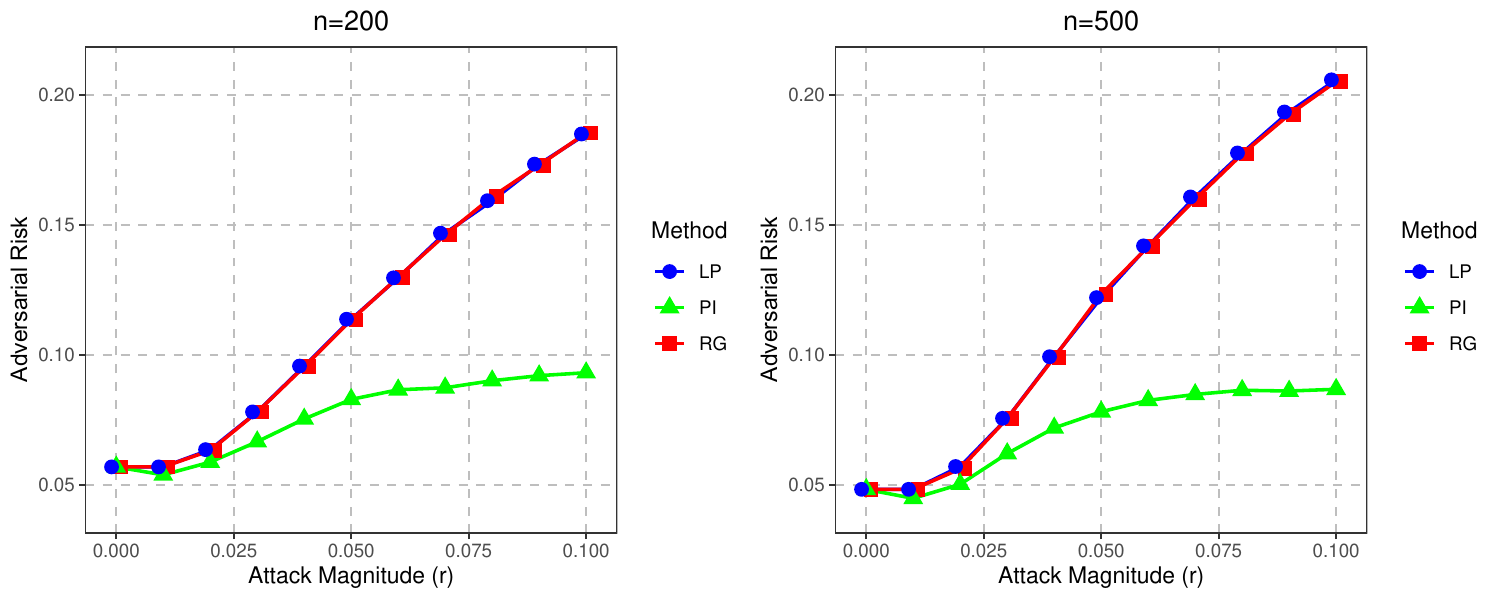}
    \end{minipage}
    }

    \subfigure[Case 2]{
    \begin{minipage}[t]{1\linewidth}
    \centering
       \includegraphics[width=5.3in]{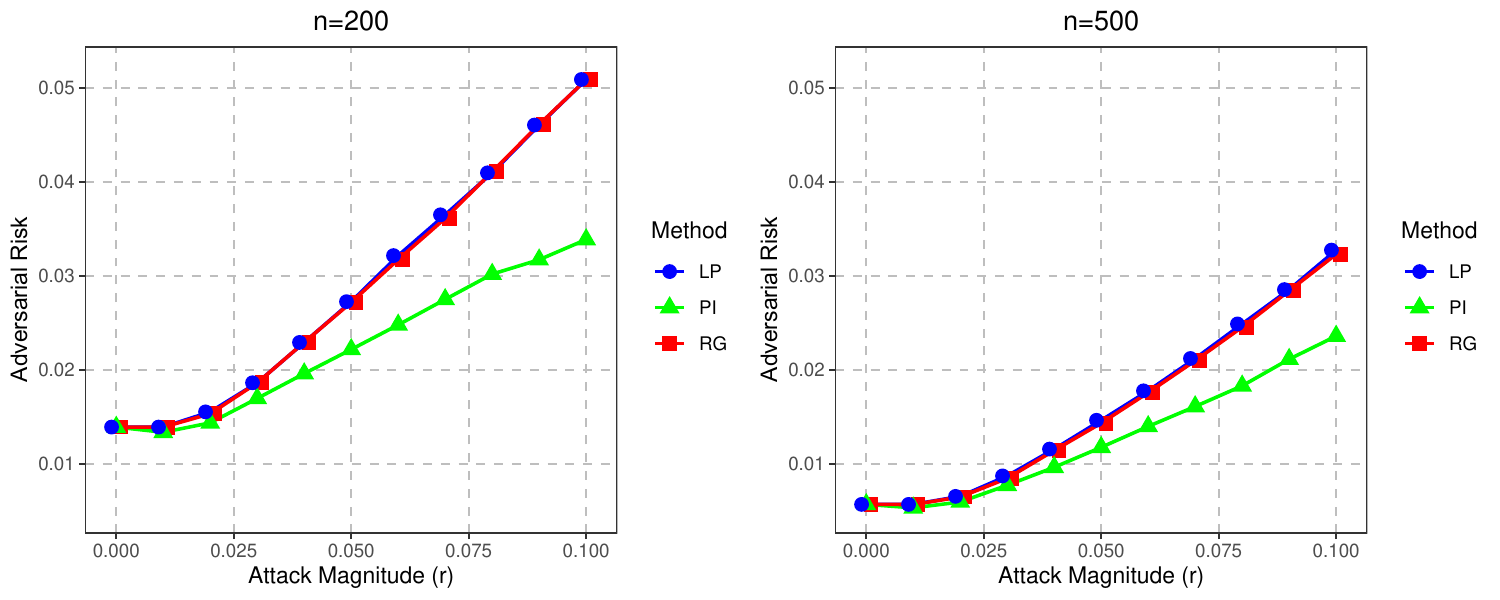}
    \end{minipage}
    }

    \caption{Adversarial risk for the three competing methods as the attack magnitude increases: panel (a) corresponds to Case 1, and panel (b) corresponds to Case 2.}
    \label{fig:s1}
\end{figure}

\begin{figure}[!htbp]
    \centering
    \subfigure[Case 3]{
    \begin{minipage}[t]{1\linewidth}
    \centering
       \includegraphics[width=5.3in]{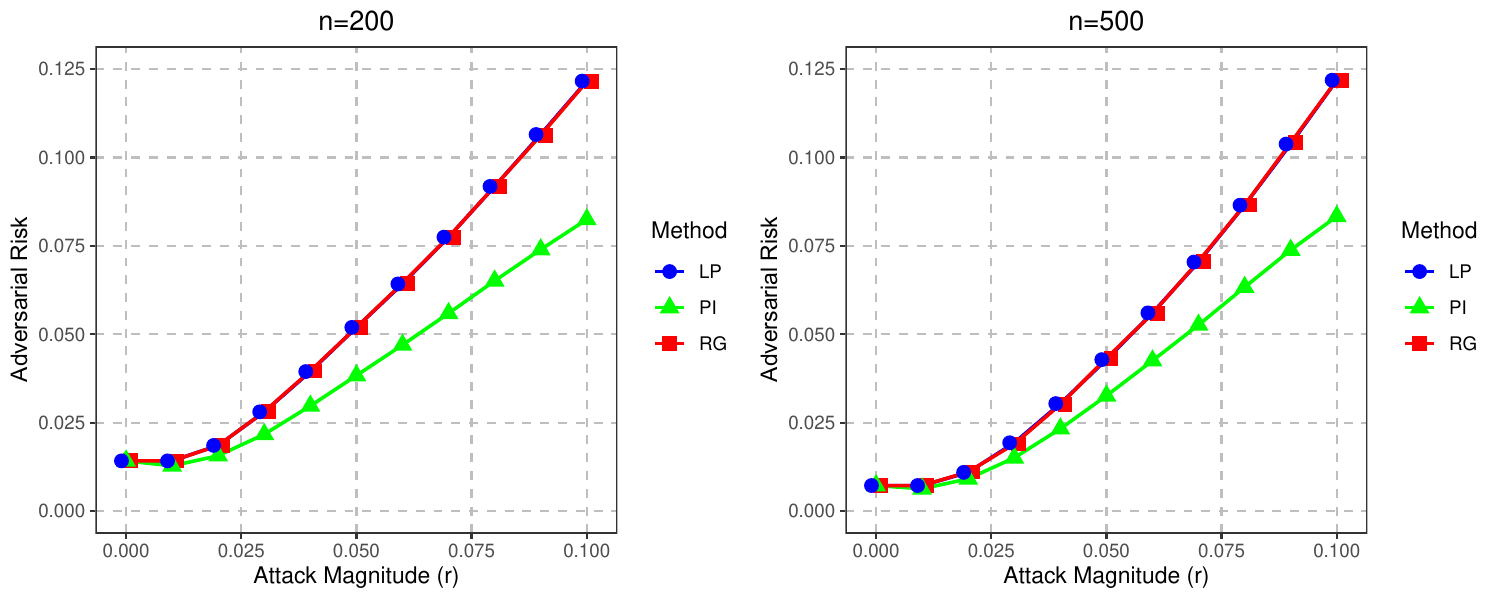}
    \end{minipage}
    }

    \subfigure[Case 4]{
    \begin{minipage}[t]{1\linewidth}
    \centering
       \includegraphics[width=5.3in]{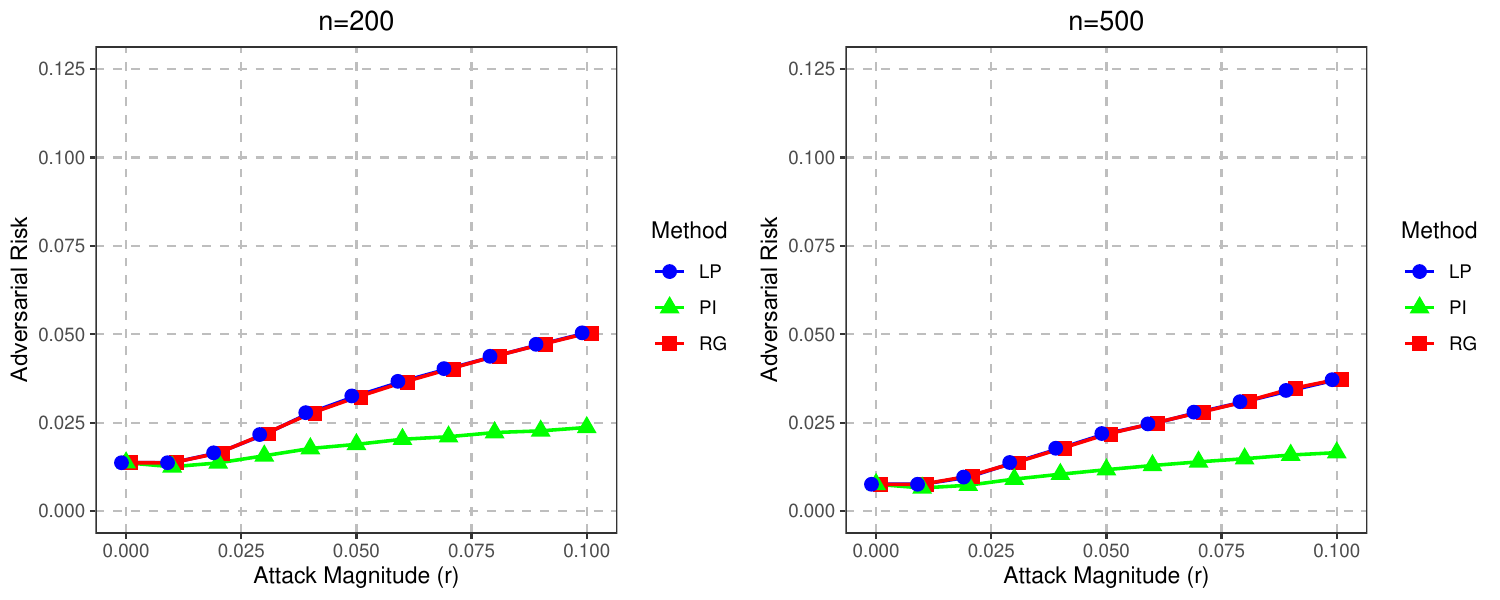}
    \end{minipage}
    }

    \caption{Adversarial risk for the three competing methods as the attack magnitude increases: panel (a) corresponds to Case 3, and panel (b) corresponds to Case 4.}
    \label{fig:s2}
\end{figure}

Building on the LP method, we consider two additional competing methods. The first (PI) follows (\ref{eq:plug_in}), where $\tilde{f}$ is the LP estimator. The second method is a ridge-type local polynomial estimator (RG), which follows the LP approach but incorporates a ridge penalty with parameter $r^2$ on the linear coefficients during the estimation of the LP coefficients. The ridge-type strategy can be seen as an approximation of adversarial training \citep{Ribeiro2023} and has also been proven to possess desirable robustness properties under several specific setups \citep{Zhang2019theoretically, Xing2021adversarially}. Figures~\ref{fig:s1}--\ref{fig:s2} present the adversarial risk for the three competing methods over 100 simulation replications. In each replication, the adversarial loss is evaluated at 100 uniformly sampled points in $[0,1]^2$.

From Figures~\ref{fig:s1}--\ref{fig:s2}, we observe a significant advantage of the PI method over the classical LP method and its ridge-type variant. For instance, in Case 2 with $n = 200$ and $r = 0.5$, the adversarial risk and its standard error for LP, RG, and PI are 2.73e-2 (0.001), 2.71e-2 (0.001), and 2.21e-2 (0.001), respectively. These results demonstrate that the adversarial plug-in procedure \ref{eq:plug_in} achieves a substantial improvement in robustness compared to the other two methods. The patterns depicted in Figures~\ref{fig:s1}--\ref{fig:s2} further corroborate the insights discussed in Sections~\ref{sec:Isotropic}--\ref{sec:anisotropic}. For example, in Case 1, where the regression function belongs to $\mathcal{F}_1(\beta,L)$ with $\beta = 1/2$, the adversarial risk curve exhibits a concave shape, consistent with the $r^{1/2}$. In Case 2, the adversarial risk curve is approximately linear as $r$ increases, which aligns with Example~\ref{exam:iso} that the adversarial loss in this case is dominated by $r$ when $r$ is large. Additionally, Figure~\ref{fig:s2} reveals that strong attacks along directions with higher variability can significantly degrade the performance of competing methods, supporting the theoretical results presented in Example~\ref{exam:aniso}.

Furthermore, although existing literature suggests that ridge-type regularization can enhance adversarial robustness under various modeling frameworks \citep[see, e.g.,][]{Zhang2019theoretically, Xing2021adversarially}, its effectiveness in the context of local nonparametric estimation remains limited. This limitation arises because ridge regularization in RG primarily controls the variation of the LP estimator at a given local point but does not regulate the variation of the estimator across different local points. Consequently, the RG method may still be vulnerable to adversarial attacks under our context.

\section{Discussion}

In this paper, we focus on the nonparametric regression problem under the adversarial attacks and examine the minimax rates of convergence in the adversarial sup-norm. Unlike the minimax analysis for the specific models in \cite{Dan2020sharp} and \cite{Xing2021adversarially}, the results established in this paper are of a general nature. They are applicable across diverse regression function classes and arbitrary additive perturbation sets. We show that the minimax rate in the adversarial setting exhibits a modular form, which equals the standard minimax rate in the absence of an adversary, plus the maximum deviation of the true function value within the perturbation set. Applying the general results to specific models is straightforward: it entails determining the standard minimax rate and calculating the largest Lipschitz constant of the functions in the target class. We further investigate two nonparametric function classes, illuminating the impacts of the different perturbation sets on the adversarial minimax rates.

It should be pointed out that the proposed adversarial plug-in estimation procedure in this paper is nonadaptive, since it depends on information about the unknown perturbation set $\Delta_n$. In the context of practical applications, an important direction for future research is to develop estimation procedures that are both adaptive across different function classes and unknown perturbation sets. Another direction is deriving the minimax rates in the general $L_p$-norm under the adversarial attacks. In the standard setting, it is well-known that the metric entropy of the regression function class plays a fundamental role in determining the minimax rates of convergence \citep{LeCam1973, birge1986estimating, yatracos1985rates, Yang1999information}. Extending these general theories to the adversarial setting is of great interest.

\appendix

\section*{Appendix}\label{appendix}
\addcontentsline{toc}{section}{Appendix}
\renewcommand{\thesection}{A.\arabic{section}}
\numberwithin{equation}{section}

\section{Proof of Lemma~\ref{lem:equ_form}}

Recall the definition of the adversarial sup-norm loss in (\ref{eq:adversarial_loss}). By the change of variable $x' = x + \delta$, we have $x' \in (x+ \Delta_n)\cap \mathcal{X}$ since $\delta \in \Delta_n$. Therefore, the adversarial loss can be expressed as
\begin{equation}\label{eq:may1}
  L_{\Delta_n}(f,\hat{f}) = \sup_{x \in \mathcal{X}} \sup_{x' \in (x+ \Delta_n)\cap \mathcal{X}} \left| f(x) - \hat{f}(x') \right|.
\end{equation}
To prove the equivalence (\ref{eq:equ_form}), it remains to show
\begin{equation}\label{eq:part_1}
  \sup_{x \in \mathcal{X}}\sup_{x' \in (x+ \Delta_n)\cap \mathcal{X}}\left|f(x)-\hat{f}(x')\right| = \sup_{x' \in \mathcal{X}}\sup_{x \in (x'-\Delta_n)\cap \mathcal{X}} \left| f(x) - \hat{f}(x') \right|.
\end{equation}

Assume that $$\sup_{x \in \mathcal{X}}\sup_{x'\in (x+ \Delta_n)\cap \mathcal{X}}|f(x)-\hat{f}(x')| > \sup_{x' \in \mathcal{X}}\sup_{x \in (x'-\Delta_n)\cap \mathcal{X}} | f(x) - \hat{f}(x') |.$$ Then there must exist $x_1 \in \mathcal{X}$ and $x_1' \in (x_1+\Delta_n)\cap \mathcal{X}$ such that $|f(x_1)-\hat{f}(x_1')| > \sup_{x' \in \mathcal{X}}\sup_{x \in (x'-\Delta_n)\cap \mathcal{X}} | f(x) - \hat{f}(x') |$. On the other hand, based on the definition of $(x_1, x_1')$, we have $x_1' \in \mathcal{X}$ and $x_1 = x_1' - \delta_1$ for some $\delta_1 \in \Delta_n$, implying that $x_1 \in (x_1'-\Delta_n)\cap \mathcal{X}$. This leads to
$$
|f(x_1)-\hat{f}(x_1')| \leq \sup_{x' \in \mathcal{X}}\sup_{x \in (x'-\Delta_n)\cap \mathcal{X}} | f(x) - \hat{f}(x') |,
$$
which is a contradiction. Likewise, we can prove that $\sup_{x \in \mathcal{X}}\sup_{x'\in (x+\Delta_n)\cap \mathcal{X}}|f(x)-\hat{f}(x')| < \sup_{x' \in \mathcal{X}}\sup_{x \in (x'-\Delta_n)\cap \mathcal{X}} | f(x) - \hat{f}(x') |$ is also impossible. Therefore, (\ref{eq:part_1}) is proved.

\section{Proof of Theorem~\ref{theo:oracle}}

Based on the results in Lemma~\ref{lem:equ_form}, we have
\begin{equation}\label{eq:L_Delta}
\begin{split}
L_{\Delta_n}(f,f') &= \sup_{x' \in \mathcal{X}}\sup_{x \in (x'-\Delta_n)\cap \mathcal{X}} \left| f(x) - f'(x') \right|. \\
\end{split}
\end{equation}
For any given $x'\in \mathcal{X}$, note that
\begin{equation}\label{eq:L_add_1}
  \begin{split}
       & \sup_{x \in (x'-\Delta_n)\cap \mathcal{X}} \left| f(x) - f'(x') \right|\\
        &= \max\left\{ \left| \sup_{x \in (x'-\Delta_n)\cap \mathcal{X}}f(x) - f'(x')\right|,  \left| \inf_{x \in (x'-\Delta_n)\cap \mathcal{X}}f(x) - f'(x')\right| \right\} \\
       & =  \left[ \left| \frac{\sup_{x \in (x'-\Delta_n)\cap \mathcal{X}}f(x) + \inf_{x \in (x'-\Delta_n)\cap \mathcal{X}}f(x)}{2} - f'(x') \right|\right.\\
      & \left. \quad+ \frac{\sup_{x \in (x'-\Delta_n)\cap \mathcal{X}}f(x) - \inf_{x \in (x'-\Delta_n)\cap \mathcal{X}}f(x)}{2} \right],
  \end{split}
\end{equation}
where the first equality follows from the fact that $|f(x) - f'(x')|$, as a piecewise linear function of $f(x)$, achieves the supremum when $f(x)$ attains either its supremum $\sup_{x \in (x'-\Delta_n)\cap \mathcal{X}}f(x)$ or its infimum $\inf_{x \in (x'-\Delta_n)\cap \mathcal{X}}f(x)$, and the second equality is established by analyzing the relative values of $\sup_{x \in (x'-\Delta_n)\cap \mathcal{X}}f(x)$, $\inf_{x \in (x'-\Delta_n)\cap \mathcal{X}}f(x)$, and $f'(x')$.

Combining (\ref{eq:L_Delta}) with (\ref{eq:L_add_1}), we obtain
\begin{equation}\label{eq:L_Delta_1}
  \begin{split}
        L_{\Delta_n}(f,f')  & = \sup_{x' \in \mathcal{X}} \left[ \left| \frac{\sup_{x \in (x'-\Delta_n)\cap \mathcal{X}}f(x) + \inf_{x \in (x'-\Delta_n)\cap \mathcal{X}}f(x)}{2} - f'(x') \right|\right.\\
      & \left.\quad + \frac{\sup_{x \in (x'-\Delta_n)\cap \mathcal{X}}f(x) - \inf_{x \in (x'-\Delta_n)\cap \mathcal{X}}f(x)}{2} \right].
  \end{split}
\end{equation}
Since $f'$ appears only in the absolute value term in (\ref{eq:L_Delta_1}), the infimum $\inf_{f'}L_{\Delta_n}(f,f')$ is therefore obtained when
\begin{equation*}
  f'(x') = f^*(x') = \frac{\sup_{x \in (x'-\Delta_n)\cap \mathcal{X}}f(x) + \inf_{x \in (x'-\Delta_n)\cap \mathcal{X}}f(x)}{2}
\end{equation*}
for any $x' \in \mathcal{X}$. And the ideal adversarial risk is given by
\begin{equation*}
  \frac{1}{2}\sup_{x' \in \mathcal{X}}\left[ \sup_{x \in (x'-\Delta_n)\cap \mathcal{X}}f(x) - \inf_{x \in (x'-\Delta_n)\cap \mathcal{X}}f(x) \right],
\end{equation*}
which completes the proof of this theorem.

\section{Proof of Theorem~\ref{theo:upper}}

From (\ref{eq:L_Delta_1}), we see
  \begin{equation}\label{eq:risk_decom}
    \begin{split}
       R_{\Delta_n}(f,\hat{f}_{\mathrm{PI}}) &= \mathbb{E}\sup_{x' \in \mathcal{X}} \left[ \left| \frac{\sup_{x \in (x'-\Delta_n)\cap \mathcal{X}}f(x) + \inf_{x \in (x'-\Delta_n)\cap \mathcal{X}}f(x)}{2} - \hat{f}_{\mathrm{PI}}(x') \right|\right.\\
      & \left. + \frac{\sup_{x \in (x'-\Delta_n)\cap \mathcal{X}}f(x) - \inf_{x \in (x'-\Delta_n)\cap \mathcal{X}}f(x)}{2} \right].
    \end{split}
  \end{equation}
  Based on the definition (\ref{eq:plug_in}) of $\hat{f}_{\mathrm{PI}}(x')$, the first term in the square bracket of (\ref{eq:risk_decom}) can be upper bounded by
  \begin{equation}\label{eq:proof1}
    \begin{split}
         & \left| \frac{\sup_{x \in (x'-\Delta_n)\cap \mathcal{X}}f(x) + \inf_{x \in (x'-\Delta_n)\cap \mathcal{X}}f(x)}{2} - \hat{f}_{\mathrm{PI}}(x') \right| \\
         & = \left| \frac{\sup_{x \in (x'-\Delta_n)\cap \mathcal{X}}f(x) + \inf_{x \in (x'-\Delta_n)\cap \mathcal{X}}f(x)}{2}\right.\\
          &\quad \left.- \frac{\sup_{x \in (x'-\Delta_n)\cap \mathcal{X}}\tilde{f}(x)+\inf_{x \in (x'-\Delta_n)\cap \mathcal{X}}\tilde{f}(x)}{2} \right|\\
         &\leq \frac{1}{2}\left| \sup_{x \in (x'-\Delta_n)\cap \mathcal{X}}f(x) - \sup_{x \in (x'-\Delta_n)\cap \mathcal{X}}\tilde{f}(x)  \right|\\
          &\quad+ \frac{1}{2}\left| \inf_{x \in (x'-\Delta_n)\cap \mathcal{X}}f(x) - \inf_{x \in (x'-\Delta_n)\cap \mathcal{X}}\tilde{f}(x)   \right|\\
         & \leq \sup_{x \in(x'-\Delta_n)\cap \mathcal{X}}\left| f(x)- \tilde{f}(x) \right| .
    \end{split}
  \end{equation}
  Combining (\ref{eq:proof1}) with (\ref{eq:risk_decom}), we have
  \begin{equation*}
    \begin{split}
         R_{\Delta_n}(f,\hat{f}_{\mathrm{PI}})& \leq \mathbb{E} \sup_{x' \in \mathcal{X}}\sup_{x \in(x'-\Delta_n)\cap \mathcal{X}}\left| f(x)- \tilde{f}(x) \right| + L_{\Delta_n}^*(f)\\
         & \leq R(f,\tilde{f}) + L_{\Delta_n}^*(f),
    \end{split}
  \end{equation*}
  where the first inequality follows from (\ref{eq:proof1})–(\ref{eq:risk_decom}) and the definition of $L_{\Delta_n}^*(f)$ in (\ref{eq:oracle_loss}), and the second inequality follows from
  \begin{equation*}
    \sup_{x' \in \mathcal{X}}\sup_{x \in(x'-\Delta_n)\cap \mathcal{X}}\left| f(x)- \tilde{f}(x) \right| \leq \sup_{x' \in \mathcal{X}}\sup_{x \in \mathcal{X}}\left| f(x)- \tilde{f}(x) \right| = \sup_{x \in \mathcal{X}}\left| f(x)- \tilde{f}(x) \right|.
  \end{equation*}
  Thus, we complete the proof of (\ref{eq:upper}).

  The second part of this theorem is proved by taking upper bound on both sides of (\ref{eq:upper}) with respect to $f \in \mathcal{F}$ and then using the condition (\ref{eq:minimax_stand}). Specifically, we have
    \begin{equation*}
    \begin{split}
       \sup_{f \in \mathcal{F}}R_{\Delta_n}(f,\hat{f}_{\mathrm{PI}}) & \lesssim \sup_{f \in \mathcal{F}}R(f,\tilde{f}) +  \sup_{f \in \mathcal{F}}L_{\Delta_n}^*(f) \\
         & \asymp \inf_{\hat{f}}\sup_{f \in \mathcal{F}}R(f,\hat{f}) +  \sup_{f \in \mathcal{F}}L_{\Delta_n}^*(f),
    \end{split}
  \end{equation*}
  which leads to (\ref{eq:minimax_upper}).

\section{Proof of Theorem~\ref{theo:lower}}

Based on the relation (\ref{eq:L_Delta_1}), we have for any $\hat{f}$,
  \begin{equation}\label{eq:lower2}
  \begin{split}
     R_{\Delta_n}(f,\hat{f}) & \geq \sup_{x' \in \mathcal{X}} \left[\frac{\sup_{x \in (x'-\Delta_n)\cap \mathcal{X}}f(x) - \inf_{x \in (x'-\Delta_n)\cap \mathcal{X}}f(x)}{2} \right] \\
       &   =   L_{\Delta_n}^*(f),
  \end{split}
  \end{equation}
  where the equality follows from (\ref{eq:oracle_loss}).
  In addition, the adversarial risk is always lower bounded by the standard risk, i.e.,
  \begin{equation}\label{eq:std2}
    R_{\Delta_n}(f,\hat{f}) = \mathbb{E}\sup_{x \in \mathcal{X}}\sup_{\delta\in \Delta_n}\left|f(x)-\hat{f}(x+\delta)\right| \geq \mathbb{E}\sup_{x \in \mathcal{X}}\left|f(x)-\hat{f}(x)\right| = R(f,\hat{f}).
  \end{equation}
  Combining (\ref{eq:lower2}) and (\ref{eq:std2}) yields the lower bound (\ref{eq:lower}). The minimax lower bound (\ref{eq:minimax_lower}) follows directly from (\ref{eq:lower}).

\section{Proof of Example~1}\label{sec:proof_exam_1}

To simplify the notation, for any $d$-dimensional multi-index $l=(l_1,l_2, \ldots, l_d) \in \mathbb{N}_0^d$, we define $|l|=l_1+l_2+\cdots+l_d$, and $l !=l_{1} ! l_{2} ! \ldots l_{d} !$. Derivatives and powers of order $l$ are denoted by $D^l=\frac{\partial^{|l|}}{\partial x^{l_1}_1\partial x^{l_2}_2\ldots\partial x^{l_d}_d}$ and $x^l=x_1^{l_1} x_2^{l_2} \ldots x_d^{l_d}$, respectively.

For any function $f$ in $\mathcal{F}_1(\beta,L)$, let
\begin{equation}\label{eq:tail_poly}
  g_k(x;t) = \sum_{|l| \leq k} \frac{D^l f(t)}{l !}(x - t)^l
\end{equation}
denote its Taylor polynomial of degree $k=\lfloor \beta \rfloor$ at point $t$. Using results from the approximation theory \citep[see, e.g.,][]{devore1993constructive}, we know that
\begin{equation}\label{eq:poly}
  \left|f(x) - g_k(x;t)\right| \leq L\sum_{|l|=k} \frac{1}{l !}|x-t|^l \cdot \left\| x - t \right\|^{\alpha},
\end{equation}
where $\alpha = \beta - k$. For completeness, we provide a simplified proof for (\ref{eq:poly}) based on the similar technique in Lemma 11.1 of \cite{gyorfi2002distribution}.
When $k = 0$, we have $\beta = \alpha$, then (\ref{eq:poly}) follows from the assumption that $f$ is $(\beta, L)$-smooth. In the case $k \geq 1$, we have
\begin{equation*}
  \begin{split}
       & \left|f(x) - g_k(x;t)\right| \\
       & = \left|f(x) - \sum_{|l| \leq k-1} \frac{D^l f(t)}{l !}(x - t)^l - \sum_{|l| = k} \frac{D^l f(t)}{l !}(x - t)^l\right|\\
       & = \left| \sum_{|l|=k} \frac{k}{l !}(x-t)^l \int_0^1(1-z)^{k-1}D^l f[t+z(x-t)] d z\right.\\
       & \left. \quad- \sum_{|l| = k} \frac{k }{l !}(x-t)^l\int_0^1(1-z)^{k-1} D^l f(t) d z \right|\\
       & = \left| \sum_{|l|=k} \frac{k}{l !}(x-t)^l \int_0^1(1-z)^{k-1}\left\{D^l f[t+z(x-t)] - D^l f(t) \right\} d z \right|\\
       & \leq L\sum_{|l|=k} \frac{1}{l !}|x-t|^l \cdot \left\| x - t \right\|^{\alpha},\\
  \end{split}
\end{equation*}
where the second equality follows from the integral form of the Taylor series remainder, and the last inequality follows from the definition of $\mathcal{F}_1(\beta,L)$.

We first construct an upper bound on $L_{\Delta_n}^*(f)$ for $f \in \mathcal{F}_1(\beta,L)$. Recall the definitions
\begin{equation*}
  \begin{split}
      2L_{\Delta_n}^*(f) &  =\sup_{x' \in \mathcal{X}}\left[ \sup_{x \in (x'-\Delta_n)\cap \mathcal{X}}f(x)-\inf_{x \in (x'-\Delta_n)\cap \mathcal{X}}f(x)\right]. \\
  \end{split}
\end{equation*}
and $r_n \triangleq \max_{\delta, \delta' \in \Delta_n }\|\delta - \delta' \|$. In addition, define $\bar{x} = (x+ x')/2$. Then we have
\begin{equation}\label{eq:key_1}
\begin{split}
2 L_{\Delta_n}^*(f) &\leq \sup_{\| x - x' \| \leq 2r_n}\left| f(x)-f(x')\right| \\
     & = \sup_{\| x - x' \| \leq 2r_n}\left| f(x) - g_k(x;\bar{x}) + g_k(x;\bar{x}) - g_k(x';\bar{x}) + g_k(x';\bar{x}) -f(x')\right|\\
     & \leq \sup_{\| x - x' \| \leq 2r_n}\left| f(x) - g_k(x;\bar{x}) \right| + \sup_{\| x - x' \| \leq 2r_n}\left|  g_k(x;\bar{x}) - g_k(x';\bar{x}) \right|\\
      &+ \sup_{\| x - x' \| \leq 2r_n}\left|  g_k(x';\bar{x}) -f(x')\right|.
\end{split}
\end{equation}
The first term at the right side of (\ref{eq:key_1}) is upper bounded by
\begin{equation}\label{eq:key_2}
\begin{split}
     \sup_{\| x - x' \| \leq 2r_n}\left| f(x) - g_k(x';\bar{x}) \right| & \leq \frac{L}{2^{\beta}}\sup_{\| x - x' \| \leq 2r_n} \sum_{|l|=k} \frac{1}{l !}\left|x-x'\right|^l \cdot \left\| x-x' \right\|^{\alpha} \\
     & \leq \frac{Lr_n^{\alpha}}{2^k}\sup_{\| x - x' \| \leq 2r_n} \sum_{|l|=k} \frac{1}{l !}\left|x-x'\right|^l \\
     & = \frac{Lr_n^{\alpha}}{2^k k!}\sup_{\| x - x' \| \leq 2r_n}\left( |x_1 - x'_1| + \cdots + |x_d - x'_d| \right)^k \\
     & \leq \frac{Lr_n^{\alpha}d^{\frac{k}{2}}}{2^k k!}\left( |x_1 - x'_1|^2 + \cdots + |x_d - x'_d|^2 \right)^\frac{k}{2} \\
     & \leq \frac{Ld^{\frac{k}{2}}r_n^{\alpha}(2r_n)^{k}}{2^{k}k!} \leq C d^{\frac{k}{2}} r_n^{\beta},
\end{split}
\end{equation}
where the first inequality follows from (\ref{eq:poly}) and the definition of $\bar{x}$, the second inequality follows from $\| x - x' \| \leq 2r_n$ and $\beta = k + \alpha$, and the third inequality follows from Jensen's inequality. The second term of (\ref{eq:key_1}) is upper bounded by
\begin{equation}\label{eq:key_3}
  \begin{split}
       &\sup_{\| x - x' \| \leq 2r_n}\left|  g_k(x;\bar{x}) - g_k(x';\bar{x}) \right|\\
       &= \sup_{\| x - x' \| \leq 2r_n}\left|  \sum_{|l| \leq k} \frac{D^l f(\bar{x})}{l !}\left[(x - \bar{x})^l - (x' - \bar{x})^l\right] \right|\\
       & = \sup_{\| x - x' \| \leq 2r_n}\left|  \sum_{|l| \leq k} \frac{D^l f(\bar{x})}{2^{|l|} l !}\left[1+(-1)^{|l|+1}\right](x - x')^l \right|\\
       & \leq C \sup_{\| x - x' \| \leq 2r_n} \sum_{1 \leq s \leq k} \sum_{|l| = s} \frac{1}{l !}|x - x'|^l\\
       & \leq  C \sup_{\| x - x' \| \leq 2r_n} \sum_{1 \leq s \leq k} \left( |x_1 - x'_1|^2 + \cdots + |x_d - x'_d|^2 \right)^\frac{s}{2}\\
       & \leq Cr_n,
  \end{split}
\end{equation}
where the first equality follows from (\ref{eq:tail_poly}), the second equality is due to the definition of $\bar{x}$, the first inequality follows from $D^l f(\bar{x})$ is bounded and $\frac{1}{2^{|l|}} \leq 1$, and the second inequality follows the similar reasoning as in the third line of (\ref{eq:key_2}). Based on the same technique in (\ref{eq:key_2}), we see the third term of (\ref{eq:key_1}) is upper bounded by
\begin{equation}\label{eq:key_4}
  \sup_{\| x - x' \| \leq 2r_n}\left|  g_k(x';\bar{x}) -f(x')\right| \leq C d^{\frac{k}{2}} r_n^{\beta}.
\end{equation}
Combining (\ref{eq:key_1}) with (\ref{eq:key_2})--(\ref{eq:key_4}), we have for all $f\in \mathcal{F}_1(\beta,L)$, $L_{\Delta_n}^*(f) \leq C d^{\frac{k}{2}} r_n^{ 1 \wedge\beta}$, i.e.,
\begin{equation*}
  \sup_{f \in \mathcal{F}_1(\beta,L)}L_{\Delta_n}^*(f) \leq C d^{\frac{k}{2}} r_n^{ 1 \wedge\beta}.
\end{equation*}

To lower bound $\sup_{f \in \mathcal{F}_1(\beta,L)}L_{\Delta_n}^*(f)$, it suffices to construct specific functions in $\mathcal{F}_1(\beta,L)$ such that $L_{\Delta_n}^*(f) \geq C r_n^{ 1 \wedge\beta}$. Given that $\Delta_n$ is a closed set, let $\delta$ and $\delta'$ be two points in $\Delta_n$ such that $\| \delta - \delta' \| = r_n$. In addition, define $D_n = \{ t\delta + (1-t)\delta': 0 \leq t \leq 1 \}$. Since there exists a $D_n$ such that $D_n \subseteq \Delta_n$, hence $L_{\Delta_n}^*(f) \geq L_{D_n}^*(f)$. Without loss of generality, we assume $D_n \subseteq \{x:x_2=x_3=\cdots =x_d = 0 \}$ and $(\delta+\delta')/2 = (1/2,0,\ldots,0)^{\top}$. Otherwise, we can construct new functions from the functions $f_1$ and $f_2$ defined below by rotations of axes and shifts of origin. Note that the rotation and transformation of a function do not change the smoothness properties of the original function. When $\beta\geq 1$, define
\begin{equation*}
  f_1(x) = L\exp(x_1-1), \quad x \in [0,1]^d.
\end{equation*}
Note that $f_1(x)$ is an infinitely differentiable function, and
\begin{equation*}
  \begin{split}
       & \left| \frac{\partial^k f_1}{\partial x_1^{k_1}\cdots \partial x_d^{k_d}}(x) - \frac{\partial^k f_1}{\partial x_1^{k_1}\cdots \partial x_d^{k_d}}(z) \right| \\
       & = \left| L\exp(x_1-1) - L\exp(z_1-1) \right|\\
       & \leq L \left| x_1- z_1\right|\leq L\left\| x - z \right\| \leq L\left\| x - z \right\|^\alpha,
  \end{split}
\end{equation*}
which verifies the conditions of $\mathcal{F}_1(\beta,L)$. Thus, $\sup_{f \in \mathcal{F}_1(\beta,L)}L_{\Delta_n}^*(f)$ is lower bounded by
\begin{equation*}
\begin{split}
     \sup_{f \in \mathcal{F}_1(\beta,L)}L_{\Delta_n}^*(f) \geq L_{D_n}^*(f_1)&\geq C\left[ \sup_{x \in \{(1/2,0,\ldots,0)^{\top}-D_n\}}f_1(x)\right.\\
     &\qquad\qquad\qquad\left.-\inf_{x \in \{(1/2,0,\ldots,0)^{\top}-D_n\}}f_1(x)\right] \\
     & \geq C r_n.
\end{split}
\end{equation*}
When $0<\beta<1$, consider the function $f_2(x) = x_1^{\beta}$. We have
\begin{equation*}
  \begin{split}
       \left| f_2(x) - f_2(z) \right| &= \left| x_1^{\beta} - z_1^{\beta} \right| \\
       & \leq \left| x_1 - z_1 \right|^{\beta} \leq \| x- z \|^{\beta}.
  \end{split}
\end{equation*}
Thus, $f_2$ belong the function class $\mathcal{F}_1(\beta,L)$ with $0<\beta <1$.
In this case, we obtain
$$
\sup_{f \in \mathcal{F}_1(\beta,L)}L_{\Delta_n}^*(f) \geq L_{D_n}^*(f_2)  \geq C r_n^{\beta},
$$
which completes the proof of this example.

\section{Proof of Example~2}
\setcounter{equation}{0}

Combining the results in \cite{bertin2004asymptotically, Bertin2004minimax}, we can obtain
\begin{equation}\label{eq:aniso_pr_1}
  \inf_{\hat{f}}\sup_{f \in \mathcal{F}_2(\beta,L)} R(f,\hat{f}) \asymp \left(\frac{\log n}{n}\right)^{\frac{\bar{\beta}}{2\bar{\beta}+d}},
\end{equation}
where $\bar{\beta}=d/(\sum_{i=1}^{d}1/\beta_i)$. Therefore, it remains to determine the rate of $\sup_{f \in \mathcal{F}_2(\beta,L)}L_{\Delta_n}^*(f)$. We first construct an upper bound on $\sup_{f \in \mathcal{F}_2(\beta,L)}L_{\Delta_n}^*(f)$. For any function $f$ in $\mathcal{F}_2(\beta,L)$, we have
\begin{equation}\label{eq:aniso_pr_2}
  \begin{split}
     2L_{\Delta_n}^*(f) & = \sup_{x' \in \mathcal{X}}\left[ \sup_{x \in (x'-\Delta_n)\cap \mathcal{X}}f(x)-\inf_{x \in (x'-\Delta_n)\cap \mathcal{X}}f(x)\right]\\
       & \leq \sup_{x' \in \mathcal{X}}\left[ \sup_{x,z\in x' - \Delta_n}\left| f(x) - f(z) \right| \right]\\
       & \leq \sup_{x' \in \mathcal{X}}\left[ \sup_{x,z\in x' - \Delta_n} \left( L_1\left| x_1 - z_1 \right|^{\beta_1} + \cdots + L_d\left| x_d - z_d \right|^{\beta_d} \right)\right]\\
       & \leq L_1r_1^{\beta_1} + \cdots + L_dr_d^{\beta_d}\\
       & \lesssim \max\{r_1^{\beta_1},\ldots, r_d^{\beta_d}\},
  \end{split}
\end{equation}
where the third step follows from the definition of $\mathcal{F}_2(\beta,L)$. Now we derive a lower bound on $\sup_{f \in \mathcal{F}_2(\beta,L)}L_{\Delta_n}^*(f)$. We just need to construct a specific function in $\mathcal{F}_2(\beta,L)$. Define $j\triangleq \arg\max_{i \in \{1,\ldots,d \}}r_i^{\beta_i}$ and a function $f_3(x)=L_{j}x_j^{\beta_{j}}$. Obviously, we have
\begin{equation*}
  \left| f_3(x) - f_3(z)\right| = L_j\left| x_{j}^{\beta_{j}} - z_{j}^{\beta_{j}} \right| \leq L_j\left| x_{j} - z_{j} \right|^{\beta_{j}}.
\end{equation*}
Thus, we see $f_3 \in \mathcal{F}_2(\beta,L)$. And $\sup_{f \in \mathcal{F}_2(\beta,L)}L_{\Delta_n}^*(f)$ is lower bounded by
\begin{equation}\label{eq:aniso_pr_3}
  \sup_{f \in \mathcal{F}_2(\beta,L)}L_{\Delta_n}^*(f) \geq L_{\Delta_n}^*(f_3) \geq L_jr_j^{\beta_{j}} \asymp \max\{r_1^{\beta_1},\ldots, r_d^{\beta_d}\}.
\end{equation}

Combining (\ref{eq:aniso_pr_1})--(\ref{eq:aniso_pr_3}) with (\ref{eq:minimax}), we have
\begin{equation*}
    \inf_{\hat{f}}\sup_{f \in \mathcal{F}_2(\beta,L)}R_{\Delta_n}(f,\hat{f}) \asymp \left(\frac{\log n}{n}\right)^{\frac{\bar{\beta}}{2\bar{\beta}+d}} + \max\{r_1^{\beta_1},\ldots, r_d^{\beta_d}\},
  \end{equation*}
  which proves the result in Example~\ref{exam:aniso}.

\newpage
\bibliographystyle{apalike}
\bibliography{mybibfile}

\end{sloppypar}
\end{document}